\documentstyle[12pt]{article}
\begin{document}
\title { On the stability of the orthogonal Pexiderized Cauchy equation }
\author{{\bf Mohammad Sal Moslehian} \\ Dept. of Math., Ferdowsi Univ.\\ P. O. Box 1159, Mashhad 91775\\ Iran\\ E-mail: msalm@math.um.ac.ir}
\date{}
\maketitle
\begin{abstract}
We investigate the stability of Pexiderized mappings in Banach modules over a unital Banach algebra. As a consequence, we establish the Hyers--Ulam stability of the orthogonal Cauchy functional equation of Pexider type $f_1(x+y)=f_2(x)+f_3(y),~~ x\perp y$ in which $\perp$ is the orthogonality in the sense of R\" atz.
\end{abstract}
\newpage

\section{Introduction.}

Assume that X is a real inner product space and $f:X\to {\bf  R}$ is a solution of the orthogonal Cauchy functional equation $f(x+y)=f(x)+f(y), <x,y>=0$. By the Pythagorean theorem $f(x)=\|x\|^2$ is a solution of the conditional equation. Of course, this function does not satisfy the additivity equation everywhere. Thus orthogonal Cauchy equation is not equivalent to the classic Cauchy equation on the whole inner product space. This phenomenon may show the significance of study of orthogonal Cauchy equation.

G. Pinsker characterized orthogonally additive functionals on an inner product space when the orthogonality is the ordinary one in such spaces $\cite{PIN}$. K. Sundaresan generalized this result to arbitrary Banach spaces equipped with the Birkhoff-James orthogonality $\cite{SUN}$. The orthogonal Cauchy functional equation $$f(x+y)=f(x)+f(y), x\perp y~~~~~(\heartsuit)$$ in which $\perp$ is an abstract orthogonality relation was first investigated by S. Gudder and D. Strawther $\cite{GU-S}$. They defined $\perp$ by a system consisting of five axioms and described the general semi-continuous real-valued solution of conditional Cauchy functional equation. In 1985, J. R\" atz introduced a new definition of orthogonality by using more restrictive axioms than of S. Gudder and D. Strawther. Moreover, he investigated the structure of orthogonally additive mappings $\cite{RAT}$. In the next step, J. R\" atz and Gy. Szab\' o investigated the problem in a rather more general framework $\cite{R-S}$.

In the recent decades, stability of functional equations have been investigated by many mathematicians. They have so many applications in information theory, Physics, Economic Theory and Social and Behaviour Sciences; cf. $\cite{ACZ}$ and $\cite{P-R}$.

The first author treating the stability of the Cauchy equation was D. H. Hyers $\cite{HYE}$ by proving that if $f$ is a mapping from a normed space $X$ into a Banach space satisfying $\|f(x+y)-f(x)-f(y)\|\leq \epsilon$ for some $\epsilon>0$, then there is a unique additive mapping $g:X\to Y$ such that $\|f(x)-g(x)\|\leq\epsilon$. Since then, the stability problem of the Cauchy equation has been extensively investigated by many mathematicians; cf. $\cite{H-I-R}$. A generalized version of Cauchy equation is the equation of Pexider type $f_1(x+y)=f_2(x)+f_3(y)$. Y.H. Lee, K.W. Jun, D.S. Shin and B.D. Kim obtained the Hyers–-Ulam–-Rassias stability of this Pexider equation; cf. $\cite{L-J}$ and $\cite{J-S-K}$. In addition, the stability of the linear and quadratic mappings in Banach modules were studied by C.-G. Park $\cite{PAR1}, \cite{PAR2}$.

R. Ger and J. Sikorska $\cite{G-S}$ investigated the orthogonal stability of the Cauchy functional equation $f(x+y)=f(x)+f(y)$, namely, they showed that if $f$ is a function from an orthogonality space $X$ into a real Banach space $Y$, $\epsilon>0$ is given and for all $x,y\in X$ with $x\perp y, f(x+y)=f(x)+f(y)$, then there exists exactly one orthogonally additive mapping $g:X\to Y$ such that for all $x\in X, \|f(x)-g(x)\|\leq \frac{16}{3}\epsilon$. 

One of the significant conditional equations is the so-called orthogonal Cauchy functional equation of Pexider type $$f_1(x+y)=f_2(x)+f_3(y),~~ x\perp y~~~~~(\diamondsuit)$$. In the present paper, we investigate the stability of Pexiderized mappings in Banach modules over a unital Banach algebra and as a consequence we establish the stability of orthogonal Pexiderized Cauchy functional equation in the spirit of Hyers--Ulam. Thus we generalize the main theorem of $\cite{G-S}$.

Throughout the paper, ${\bf R}$ and ${\bf R_+}$ denote the sets of real and nonnegative real numbers, respectively. $A$ is a unital real Banach algebra with unit $1$ and unit sphere $A_1$. In addition, all modules are assumed to be unit linked real left modules over $A$. The reader is referred to $\cite{B-D}$ for more details on the theory of Banach modules.

\section {Preliminaries.}

There are several orthogonality notions on a real normed space such as Birkhoff-James, Boussouis, (semi-)inner product, Singer, Carlsson, area, unitary-Boussouis, Roberts, Phythagorean, isoscelesa and Diminnie (see e.g. $\cite{A-B1}$ and $\cite{A-B2}$).

Let us recall the orthogonality in the sense of J. R\" atz; cf. $\cite{RAT}$. 

Suppose $X$ is a real vector space (algebraic module) with $\dim X\geq 2$ and $\perp$ is a binary relation on $X$ with the following properties:\\
(O1) totality of $\perp$ for zero: $x\perp 0, 0\perp x$ for all $x\in X$;\\
(O2) independence: if $x,y\in X-\{0\}, x\perp y$, then $x,y$ are linearly independent;\\
(O3) homogeneity: if $x,y\in X, x\perp y$, then $\alpha x\perp \beta y$ for all $\alpha, \beta\in {\bf R}$;\\
(O4) the Thalesian property: if $P$ is a $2$-dimensional subspace of $X, x\in P$ and $\lambda\in {\bf R_+}$, then there exists $y_0\in P$ such that $x\perp y_0$ and $x+y_0\perp\lambda x-y_0$.\\
The pair $(X,\perp)$ is called an orthogonality space (module). By an orthogonality normed space (normed module) we mean an orthogonality space (module) having a normed (normed module) structure.

Some interesting examples are (i) The trivial orthogonality on a vector space $X$ defined by (O1), and for non-zero elements $x,y\in X$, $x\perp y$ if and only if $x,y$ are linearly independent.\\
(ii) The ordinary orthogonality on an inner product space $(X, \langle.,.\rangle)$ given by $x\perp y$ if and only if $\langle x,y\rangle=0$.\\
(iii) The Birkhoff-James orthogonality on a normed space $(X,\|.\|)$ defined by $x\perp y$ if and only if $\|x+\lambda y\|\geq \|x\|$ for all $\lambda\in {\bf R}$.

Let $X$ be a vector space (an orthogonality space) and $(Y,+)$ be an abelian group. A mapping $f:X\to Y$ is called (orthogonally) additive if it satisfies the so-called (orthogonal) Cauchy function equation $f(x+y)=f(x)+f(y)$ for all $x,y\in X$ (with $x\perp y$). Further, if $X$ and $Y$ are modules and $f(ax)=af(x)$ for all $a\in A$ and $x\in X$, then $f$ is called $A$-linear. 
A mapping $f:X\to Y$ is said to be (orthogonally) quadratic if it satisfies the so-called (orthogonally) quadratic function equation $f(x+y)+f(x-y)=2f(x)+2f(y)$ for all $x,y\in X$ (with $x\perp y$). Further, if $X$ and $Y$ are modules and $f(ax)=a^2f(x)$ for all $a\in A$ and $x\in X$, then $f$ is called $A$-quadratic. 

In 1985, R\" atz gave the following significant result (cf. Corollary 7 of $\cite{RAT}$).

{\bf Theorem (*).} If $(Y,+)$ is uniquely $2$-devisable (i.e. the mapping $\omega:Y\to Y, \omega(y)=2y$ is bijective), in particular a vector space, then every solution $f$ of the orthogonally additive function $(\heartsuit)$ equation has the form $f=Q+T$ with $Q$ quadratic and $T$ additive. 

\section {Orthogonal Stability in Banach modules.}

In this section, applying some ideas from $\cite{G-S}, \cite{J-S}, \cite{PAR1}$ and using sequences of Hyers' type $\cite{HYE}$ being a useful tool in the theory of stability of equations, among several things, we deal with the conditional stability problem for equation $(\diamondsuit)$.

{\bf Lemma 1.} Suppose $(X,\perp)$ is an orthogonality module and $(Y, \|.\|)$ is a real Banach module. Let $F_1,F_2, F_3:X\to Y$ be even mappings fulfilling 
\begin{equation}
\|F_1(ax+ay)-a^2F_2(x)-a^2F_3(y)\|\leq\epsilon
\end{equation}
for some $\epsilon$, for all $a\in A_1$ and for all $x,y\in X$ with $x\perp y$. Assume that $F_i(0)=0, i=1,2,3$. Then there exists a unique quadratic mapping $Q:X\to Y$ such that
\begin{eqnarray*}
\|F_1(x)-Q(x)\|\leq\frac{13}{3}\epsilon\\
\|F_2(x)-Q(x)\|\leq\frac{16}{3}\epsilon\\
\|F_3(x)-Q(x)\|\leq\frac{16}{3}\epsilon
\end{eqnarray*}
for all $x\in X$. Moreover, $Q(ax)=a^2Q(x)$ for all $a\in A_1, x\in X$.

{\bf Proof.} For every $x\in X, x\perp 0$. So we can put $a=1$ and $y=0$ in $(1)$ to obtain 
\begin{equation}
\|F_1(x)-F_2(x)\|\leq \epsilon,~~x\in X
\end{equation}
Similarly, we can put $a=1$ and $x=0$ in $(1)$ to obtain 
\begin{equation}
\|F_1(y)-F_3(y)\|\leq \epsilon,~~y\in X
\end{equation}
If $x\perp y$, then by (O3) $x\perp -y$. Hence we can put $a=1$ and replace $y$ by $-y$ in $(1)$ to get
\begin{equation}
\|F_1(x-y)-F_2(x)-F_3(y)\|\leq \epsilon ~~~~~x\perp y.
\end{equation}
Let $a\in A_1$ and $x\in X$ be fixed. By (O4) there exists $y_0\in X$ such that $x\perp y_0$ and $x+y_0\perp x-y_0$. Replacing $x$ and $y$ by $x+y_0$ and $x-y_0$ in $(1)$, we have
\begin{equation}
\|F_1(2ax)-a^2F_2(x+y_0)-a^2F_3(x-y_0)\|\leq \epsilon.
\end{equation}
By (O3), $\frac{x+y_0}{2}\perp\pm\frac{x-y_0}{2}$ and so by using $(1)$ with $a=1$, we obtain
\begin{eqnarray*}
\|F_1(x)-F_2(\frac{x+y_0}{2})-F_3(\frac{x-y_0}{2})\|\leq \epsilon.
\end{eqnarray*}
\begin{eqnarray*}
\|F_1(y_0)-F_2(\frac{x+y_0}{2})-F_3(\frac{y_0-x}{2})\|\leq \epsilon
\end{eqnarray*}
whence, by virtue of triangular inequality, we get
\begin{equation}
\|F_1(y_0)-F_1(x)\|\leq 2\epsilon.
\end{equation}
It follows from 
\begin{eqnarray*}
\|F_1(2ax)-a^2F_1(x+y_0)-a^2F_1(x-y_0)\|\leq\\
\|F_1(2ax)-a^2F_2(x+y_0)-a^2F_3(x-y_0)\|+\\
\|a^2F_2(x+y_0)-a^2F_1(x+y_0)\|+\|a^2F_3(x-y_0)-a^2F_1(x-y_0)\|
\end{eqnarray*}
and $(2), (3)$ and $(5)$ that
\begin{equation}
\|F_1(2ax)-a^2F_1(x+y_0)-a^2F_1(x-y_0)\|\leq 3\epsilon.
\end{equation}
It follows from 
\begin{eqnarray*}
\|F_1(x-y_0)+F_1(x+y_0)-4F_1(x)\|\leq \|F_1(x-y_0)-F_2(x)-F_3(y_0)\|\\
+\|F_1(x+y_0)-F_2(x)-F_3(y_0)\|+2\|F_2(x)-F_1(x)\|\\+2\|F_3(y_0)-F_1(y_0)\|+2\|F_1(y_0)-F_1(x)\|
\end{eqnarray*}
and $(1), (2), (3), (4)$ and $(6)$ that
\begin{equation}
\|F_1(x-y_0)+F_1(x+y_0)-4F_1(x)\|\leq 10\epsilon.
\end{equation}
It follows from $(7),(8)$ and
\begin{eqnarray*}
\|F_1(2ax)-4a^2F_1(x)\|&\leq& \|F_1(2ax)-a^2F_1(x+y_0)-a^2F_1(x-y_0)\|\\
&+&\|a^2F_1(x-y_0)+a^2F_1(x+y_0)-4a^2F_1(x)\|
\end{eqnarray*}
that
\begin{equation}
\|F_1(2ax)-4a^2F_1(x)\|\leq 13\epsilon.
\end{equation}
Putting $a=1$ in $(9)$ and using induction we infer that
\begin{equation}
\|4^{-n}F_1(2^nx)-F_1(x)\|\leq (1-\frac{1}{4^n})\frac{13\epsilon}{3}.
\end{equation}
Hence $\{4^{-n}F_1(2^nx)\}$ is a Cauchy sequence in the Banach space $Y$ and so is convergent. Set $\phi(x):=\displaystyle{\lim_{n\to\infty}}4^{-n}F_1(2^nx)$. By $(10)$, $\|\phi(x)-F_1(x)\|\leq \frac{13\epsilon}{3}$. Applying inequality $(2)$, we get $\|4^{-n}F_1(2^nx)-4^{-n}F_2(2^nx)\|\leq\frac{\epsilon}{4^n}$ whence $\phi(x)=\displaystyle{\lim_{n\to\infty}}4^{-n}F_2(2^nx)$. Similarly, it follows from $(3)$ that $\phi(x)=\displaystyle{\lim_{n\to\infty}}4^{-n}F_3(2^nx)$.

For all $x,y\in X$ with $x\perp y$, inequality $(1)$ yields 
$$\|4^{-n}F_1(2^n(x+y))-4^{-n}F_2(2^nx)-4^{-n}F_3(2^ny)\|\leq 4^{-n}\epsilon.$$

Taking the limit, we deduce that $\phi(x+y)-\phi(x)-\phi(y)=0$. Hence $\phi$ is orthogonally additive. Theorem (*) states that $\phi$ can be expressed as the sum $Q+S$ of two quadratic and additive mappings. Hence $\|Q(x)+S(x)-F_1(x)\|\leq \frac{13\epsilon}{3}$. Since $F_1$ is an even function and $Q(-x)=Q(x)$, we have $\|S(x)\|\leq\frac{1}{2}\|Q(x)+S(x)-F_1(x)\|+\frac{1}{2}\|-Q(-x)-S(-x)+F_1(-x)\|\leq \frac{13\epsilon}{3}$. Thus $\|Sx\|=\frac{1}{n}\|S(nx)\|\leq\frac{13\epsilon}{3n}$ for all $n$.Therefore $Sx=0$ and so $\phi(x)=Q(x)$. Thus
\begin{equation}
\|Q(x)-F_1(x)\|\leq \frac{13\epsilon}{3}.
\end{equation}
Moreover, the inequality $(9)$ yields
\begin{eqnarray*}
\|F_1(2^nax)-4a^2F_1(2^{n-1}x)\|\leq 13\epsilon
\end{eqnarray*}
for all $x\in X, a\in A_1$ and so 
\begin{eqnarray*}
Q(ax)=\displaystyle{\lim_{n\to\infty}}4^{-n}F_1(2^nax)=\displaystyle{\lim_{n\to\infty}}4^{-(n-1)}a^2F_1(2^{n-1}x)=a^2Q(x).
\end{eqnarray*}

If $Q':X\to Y$ is another quadratic mapping fulfilling $\|Q'(x)-F_1(x)\|\leq \frac{13\epsilon}{3}$, then $\|Q(x)-Q'(x)\|\leq\frac{1}{n^2}(\|Q(nx)-F_1(nx)\|+\|Q'(nx)-F_1(nx)\|)\leq \frac{26\epsilon}{3n^2}$. Tending $n$ to $\infty$ we get $Q=Q'$ which proves the uniqueness assertion.
Further, inequalities $(2), (3), (11)$ imply that
\begin{eqnarray*}
\|F_2(x)-Q(x)\|\leq \|F_2(x)-F_1(x)\|+\|F_1(x)-Q(x)\|\leq\epsilon+\frac{13\epsilon}{3}=\frac{16\epsilon}{3}
\end{eqnarray*}
and
\begin{eqnarray*}
\|F_3(x)-Q(x)\|\leq \|F_3(x)-F_1(x)\|+\|F_1(x)-Q(x)\|\leq\epsilon+\frac{13\epsilon}{3}=\frac{16\epsilon}{3}.\Box
\end{eqnarray*}

{\bf Remark 1.} In the proof of Lemma 1 we do not use the assumptions that $F_2$ is even and $F_1(0)=0$.

{\bf Corollary 1.} Suppose $(X,\perp)$ is an orthogonality space and $(Y, \|.\|)$ is a real Banach space. Let $F_1,F_2, F_3:X\to Y$ be even mappings fulfilling 
\begin{eqnarray*}
\|F_1(x+y)-F_2(x)-F_3(y)\|\leq\epsilon
\end{eqnarray*}
for some $\epsilon$ and for all $x,y\in X$ with $x\perp y$. Assume that $F_i(0)=0, i=1,2,3$. Then there exists a unique quadratic mapping $Q:X\to Y$ such that
\begin{eqnarray*}
\|F_1(x)-Q(x)\|\leq\frac{13}{3}\epsilon\\
\|F_2(x)-Q(x)\|\leq\frac{16}{3}\epsilon\\
\|F_3(x)-Q(x)\|\leq\frac{16}{3}\epsilon
\end{eqnarray*}
for all $x\in X$. 

{\bf Proof.} Use the same reasoning as in the proof of Lemma 1 with $a=1$.$\Box$

{\bf Lemma 2.} Suppose $(X,\perp)$ is an orthogonality module and $(Y, \|.\|)$ is a real Banach module. Let $F_1,F_2, F_3:X\to Y$ be odd mappings fulfilling 
\begin{equation}
\|F_1(ax+ay)-aF_2(x)-aF_3(y)\|\leq\epsilon
\end{equation}
for some $\epsilon$, for all $a\in A_1$ and for all $x,y\in X$ with $x\perp y$. Then there exists a unique additive mapping $T:X\to Y$ such that
\begin{eqnarray*}
\|F_1(x)-T(x)\|\leq 7\epsilon\\
\|F_2(x)-T(x)\|\leq 8\epsilon\\
\|F_3(x)-T(x)\|\leq 8\epsilon
\end{eqnarray*}
for all $x\in X$. Moreover, $T(ax)=aT(x)$ for all $a\in A_1, x\in X$.

{\bf Proof.} For every $x\in X, x\perp 0$. So we can put $a=1$ and $y=0$ in $(12)$ to obtain 
\begin{equation}
\|F_1(x)-F_2(x)\|\leq \epsilon,~~x\in X
\end{equation}
Similarly we can put $a=1$ and $x=0$ in $(12)$ to obtain 
\begin{equation}
\|F_1(y)-F_3(y)\|\leq \epsilon,~~y\in X
\end{equation}
If $x\perp y$, then by (O3) $x\perp -y$. Hence we can put $a=1$ and replace $y$ by $-y$ in $(12)$ to get
\begin{equation}
\|F_1(x-y)-F_2(x)+F_3(y)\|\leq \epsilon~~~~~x\perp y.
\end{equation}
Let $a\in A_1$ and $x\in X$ be fixed. By (O4) there exists $y_0\in X$ such that $x\perp y_0$ and $x+y_0\perp x-y_0$. It follows from $(12)$ that
\begin{equation}
\|F_1(2ax)-aF_2(x+y_0)-aF_3(x-y_0)\|\leq \epsilon.
\end{equation}
It follows from 
\begin{eqnarray*}
\|F_1(2ax)-aF_1(x+y_0)-aF_1(x-y_0)\|\leq \|F_1(2ax)-aF_2(x+y_0)-aF_3(x-y_0)\|\\
+\|aF_2(x+y_0)-aF_1(x+y_0)\|+\|aF_3(x-y_0)-aF_1(x-y_0)\|
\end{eqnarray*}
and $(13), (14)$ and $(16)$ that
\begin{equation}
\|F_1(2ax)-aF_1(x+y_0)-aF_1(x-y_0)\|\leq 3\epsilon.
\end{equation}
It follows from 
\begin{eqnarray*}
\|F_1(x+y_0)+F_1(x-y_0)-2F_1(x)\|\leq \|F_1(x+y_0)-F_2(x)-F_3(y_0)\|\\
+\|F_1(x-y_0)-F_2(x)+F_3(y_0)\|+2\|F_2(x)-F_1(x)\|
\end{eqnarray*}
and $(12), (13)$ and $(15)$ that
\begin{equation}
\|F_1(x+y_0)+F_1(x-y_0)-2F_1(x)\|\leq 4\epsilon.
\end{equation}
Now $(17)$ and $(18)$ and 
\begin{eqnarray*}
\|F_1(2ax)-2aF_1(x)\|&\leq&\|F_1(2ax)-aF_1(x+y_0)-aF_1(x-y_0)\|\\
&+&\|aF_1(x+y_0)+aF_1(x-y_0)-2aF_1(x)\|
\end{eqnarray*}
yield
\begin{equation}
\|F_1(2ax)-2aF_1(x)\|\leq 7\epsilon.
\end{equation}
Putting $a=1$ in $(19)$ and using induction we infer that
\begin{equation}
\|2^{-n}F_1(2^nx)-F_1(x)\|\leq (1-\frac{1}{2^n})7\epsilon.
\end{equation}
Hence $\{2^{-n}F_1(2^nx)\}$ is a Cauchy sequence in the Banach space $Y$ and so is convergent. Set $\psi(x):=\displaystyle{\lim_{n\to\infty}}2^{-n}F_1(2^nx)$. By $(20)$, $\|\psi(x)-F_1(x)\|\leq 7\epsilon$. Applying inequality $(13)$, we get $\|2^{-n}F_1(2^nx)-2^{-n}F_2(2^nx)\|\leq\frac{\epsilon}{2^n}$ whence $\psi(x)=\displaystyle{\lim_{n\to\infty}}2^{-n}F_2(2^nx)$. Similarly, it follows from $(14)$ that $\psi(x)=\displaystyle{\lim_{n\to\infty}}2^{-n}F_3(2^nx)$.

For all $x,y\in X$ with $x\perp y$, inequality $(12)$ yields 
$$\|2^{-n}F_1(2^n(x+y))-2^{-n}F_2(2^nx)-2^{-n}F_3(2^ny)\|\leq 2^{-n}.\epsilon.$$

Taking the limit, we deduce that $\psi(x+y)-\psi(x)-\psi(y)=0$. Hence $a$ is orthogonally additive. Theorem (*) states that $a$ can be expressed as the sum $P+T$ of two quadratic and additive mappings. Hence $\|P(x)+T(x)-F_1(x)\|\leq 7\epsilon$. Since $F_1$ is an odd function and $T(-x)=-T(x)$, we have $\|P(x)\|\leq\frac{1}{2}\|P(x)+T(x)-F_1(x)\|+\frac{1}{2}\|P(-x)+T(-x)-F_1(-x)\|\leq 7\epsilon$. Thus $\|Px\|=\frac{1}{n^2}\|P(nx)\|\leq\frac{7\epsilon}{n^2}$ for all $n$.Therefore $Px=0$ and so $\psi(x)=T(x)$. Thus
\begin{equation}
\|T(x)-F_1(x)\|\leq 7\epsilon
\end{equation}

Moreover, the inequality $(19)$ yields
\begin{eqnarray*}
\|F_1(2^nax)-2aF_1(2^{n-1}x)\|\leq 7\epsilon
\end{eqnarray*}
for all $x\in X, a\in A_1$ and so 
\begin{eqnarray*}
T(ax)=\displaystyle{\lim_{n\to\infty}}2^{-n}F_1(2^nax)=\displaystyle{\lim_{n\to\infty}}2^{-(n-1)}aF_1(2^{n-1}x)=aT(x).
\end{eqnarray*}

If $T':X\to Y$ is another additive mapping fulfilling $\|T'(x)-F_1(x)\|\leq 7\epsilon$, then $\|T(x)-T'(x)\|\leq\frac{1}{n}(\|T(nx)-F_1(nx)\|+\|T'(nx)-F_1(nx)\|)\leq \frac{14\epsilon}{n}$. Tending $n$ to $\infty$ we infer that $T=T'$ which proves the uniqueness assertion.
Further, inequalities $(13), (14), (21)$ imply that
\begin{eqnarray*}
\|F_2(x)-T(x)\|\leq \|F_2(x)-F_1(x)\|+\|F_1(x)-T(x)\|\leq\epsilon+7\epsilon=8\epsilon
\end{eqnarray*}
and
\begin{eqnarray*}
\|F_3(x)-T(x)\|\leq \|F_3(x)-F_1(x)\|+\|F_1(x)-T(x)\|\leq\epsilon+7\epsilon=8\epsilon.\Box
\end{eqnarray*}

{\bf Remark 2.} In the proof of Lemma 2 we do not use the assumption that $F_2$ is odd.

{\bf Corollary 2.} Suppose $(X,\perp)$ is an orthogonality space and $(Y, \|.\|)$ is a real Banach space. Let $F_1,F_2, F_3:X\to Y$ be odd mappings fulfilling 
\begin{eqnarray*}
\|F_1(x+y)-F_2(x)-F_3(y)\|\leq\epsilon
\end{eqnarray*}
for some $\epsilon$ and for all $x,y\in X$ with $x\perp y$. Then there exists a unique additive mapping $T:X\to Y$ such that
\begin{eqnarray*}
\|F_1(x)-T(x)\|\leq 7\epsilon\\
\|F_2(x)-T(x)\|\leq 8\epsilon\\
\|F_3(x)-T(x)\|\leq 8\epsilon
\end{eqnarray*}
for all $x\in X$. 

{\bf Proof.} Use the same reasoning as in the proof of Lemma 2 with $a=1$.$\Box$.

{\bf Theorem 1.} Suppose $(X,\perp)$ is an orthogonality module and $(Y, \|.\|)$ is a real Banach module. Let $f_1, f_2, f_3:X\to Y$ be mappings fulfilling 
\begin{equation}
\|f_1(ax+ay)-abf_2(x)-abf_3(y)\|\leq\epsilon
\end{equation}
for some $\epsilon$, all $a,b\in A_1$ and for all $x,y\in X$ with $x\perp y$. Then there exists exactly a quadratic mapping $Q:X\to Y$ and an additive mapping $T:X\to Y$ such that
\begin{eqnarray*}
\|f_1(x)-f_1(0)-Q(x)-T(x)\|\leq \frac{68}{3}\epsilon\\
\|f_2(x)-f_2(0)-Q(x)-T(x)\|\leq \frac{80}{3}\epsilon\\
\|f_3(x)-f_3(0)-Q(x)-T(x)\|\leq \frac{80}{3}\epsilon
\end{eqnarray*}
for all $x\in X$. Furthermore, $T(ax)=aT(x)$ and $Q(ax)=a^2Q(x)$ for all $x\in X, a\in A_1$.

{\bf Proof.} For $1\leq i\leq 3$ define $F_i(x)=f_i(x)-f_i(0)$ and denote the even and odd parts of $F_i$ by $F_i^e, F_i^o$, respectively. Clearly $F_i^e(0)=F_i^o(0)=F_i(0)=0, i=1,2,3$

Putting $x=y=0$ in $(22)$ and subtracting the argument of the norm of the resulting inequality from that of inequality $(22)$ we get
\begin{equation}
\|F_1(ax+ay)-abF_2(x)-abF_3(y)\|\leq 2\epsilon.
\end{equation}
If $x\perp y$ then, by (O3), $-x\perp -y$. Hence we can replace $x$ by $-x$ and $y$ by $-y$ in $(23)$ to obtain
\begin{equation}
\|F_1(-ax-ay)-abF_2(-x)-abF_3(-y)\|\leq 2\epsilon.
\end{equation}
By virtue of triangular inequality and $(23)$ and $(24)$ we have
\begin{equation}
\|F_1^e(ax+ay)-abF_2^e(x)-abF_3^e(y)\|\leq 2\epsilon
\end{equation}
\begin{equation}
\|F_1^o(ax+ay)-abF_2^o(x)-abF_3^o(y)\|\leq 2\epsilon
\end{equation}
for all $x,y\in X$.\\
Putting $a=b$ in $(25)$ and applying Lemma 1, there exists a quadratic mapping $Q$ such that
$\|F_1^e(x)-Qx\|\leq 14\epsilon$ and $Q(ax)=a^2Q(x)$ for all $x\in X, a\in A_1$.

Putting $b=1$ in $(26)$ and applying Lemma 2, there exists an additive mapping $T$ such that
$\|F_1^o(x)-Tx\|\leq \frac{26}{3}\epsilon$ and $T(ax)=aT(x)$ for all $x\in X, a\in A_1$.

Hence
$$\|f_1(x)-f_1(0)-Q(x)-T(x)\|\leq \|F_1^e(x)-Qx\|+\| F_1^o(x)-Tx\|\leq 14\epsilon+\frac{26}{3}\epsilon=\frac{26}{3}\epsilon$$
Similarly, one can shows that 
$$\|f_2(x)-f_2(0)-Q(x)-T(x)\|\leq \frac{80}{3}\epsilon$$
$$\|f_3(x)-f_3(0)-Q(x)-T(x)\|\leq \frac{80}{3}\epsilon$$
Using the same method as the proof of Lemmas 1 and 2, the rest can be easily proved.$\Box$

{\bf Remark 3.} If we replace condition $(22)$ by
$$\|f_1(ax+ay)-af_2(x)-af_3(y)\|\leq\epsilon,~~~ x\perp y, a=a^2, \|a\|=1$$
then Theorem 1 is still true except that $T(ax)=aT(x)$ and $Q(ax)=a^2Q(x)$ hold merely for idempotents $a\in A_1$. This may be of special interest whenever we deal with the Banach algebras generated by their idempotents.

{\bf Remark 4.} If $f_2=\alpha f_1$ for some scalar $\alpha\neq 1$ then by $(2)$ and $(13)$ we have $|1-\alpha|\|F_1^e(x)\|\leq \epsilon$ and $|1-\alpha|\|F_1^o(x)\|\leq \epsilon$ for all $x\in X$. Hence $Q(x)=\displaystyle{\lim_{n\to\infty}}4^{-n}F_1^e(2^nx)=0$ and $T(x)=\displaystyle{\lim_{n\to\infty}}2^{-n}F_1^o(2^nx)=0$ for all $x\in X$. In particular, it follows from the conclusions of Theorem 1 that 
\begin{eqnarray*}
\|f_1(x)\|\leq \|f_1(0)\|+\frac{68}{3}\epsilon\\
\|f_3(x)\|\leq \|f_3(0)\|+\frac{80}{3}\epsilon
\end{eqnarray*}
There is a similar assertion when $f_3=\alpha f_1$ for some scalar $\alpha\neq 1$.

{\bf Theorem 2.} Suppose $(X,\perp)$ is an orthogonality normed module and $(Y, \|.\|)$ is a real Banach module. Let $f_1, f_2, f_3:X\to Y$ be mappings fulfilling 
\begin{equation}
\|f_1(ax+ay)-abf_2(x)-abf_3(y)\|\leq\epsilon
\end{equation}
for some $\epsilon$, all $a,b\in A_1$ and for all $x,y\in X$ with $x\perp y$. Then there exist exactly a quadratic mapping $Q:X\to Y$ and a additive mapping $T:X\to Y$ such that
\begin{eqnarray*}
\|f_1(x)-f_1(0)-Q(x)-T(x)\|\leq \frac{68}{3}\epsilon\\
\|f_2(x)-f_2(0)-Q(x)-T(x)\|\leq \frac{80}{3}\epsilon\\
\|f_3(x)-f_3(0)-Q(x)-T(x)\|\leq \frac{80}{3}\epsilon
\end{eqnarray*}
for all $x\in X$. In addition, if the mapping $t\mapsto f_1(tx)$ is continuous for each fixed $x\in X$, then $T$ is $A$-linear and $Q$ is $A$-quadratic.

{\bf Proof.} we use the notation of Proof of Theorem 1.

By Theorem 1, there there exist exactly a quadratic mapping $Q(x)=\displaystyle{\lim_{n\to\infty}}4^{-n}F_1^e(2^nx)$ and an additive mapping $T(x)=\displaystyle{\lim_{n\to\infty}}2^{-n}F_1^o(2^nx)$ satisfying the inquired inequalities and as well $$T(ax)=aT(x), Q(ax)=a^2Q(x), x\in X, a\in A_1$$.
For each fixed $x\in X$, because of the continuity of $t\mapsto f_1(tx)$, we deduce that $t\mapsto F_1^e(tx)$ and $t\mapsto F_1^o(tx)$ are continuous too. By the same arguing as in the proof of the theorem of $\cite{RAS}$, we can establish that $T$ is ${\bf R}$-linear and $Q$ is ${\bf R}$-quadratic.

Now for all $a\in A, x\in X$ we have
\begin{eqnarray*}
Q(ax)=Q(|a|\frac{a}{|a|}x)=|a|^2Q(\frac{a}{|a|}x)=|a|^2\frac{a^2}{|a|^2}Q(x)=a^2Q(x)
\end{eqnarray*}
and similarly, 
\begin{eqnarray*}
T(ax)=T(|a|\frac{a}{|a|}x)=|a|T(\frac{a}{|a|}x)=|a|\frac{a}{|a|}T(x)=aT(x).\Box
\end{eqnarray*}

{\bf Corollary 3.} Suppose $(X,\perp)$ is an orthogonality complex normed space and $(Y, \|.\|)$ is a complex Banach space. Let $f_1, f_2, f_3:X\to Y$ be mappings fulfilling 
\begin{equation}
\|f_1(\lambda x+\lambda y)-\lambda\mu f_2(x)-\lambda\mu f_3(y)\|\leq\epsilon
\end{equation}
for some $\epsilon$, all $\lambda,\mu\in\{z\in{\bf C}: |z|=1\}$ and for all $x,y\in X$ with $x\perp y$. If the mapping $t\mapsto f_1(tx)$ is continuous for each fixed $x\in X$, then there exist exactly a {\bf C}-quadratic mapping $Q:X\to Y$ and a {\bf C}-additive mapping $T:X\to Y$ such that
\begin{eqnarray*}
\|f_1(x)-f_1(0)-Q(x)-T(x)\|\leq \frac{68}{3}\epsilon\\
\|f_2(x)-f_2(0)-Q(x)-T(x)\|\leq \frac{80}{3}\epsilon\\
\|f_3(x)-f_3(0)-Q(x)-T(x)\|\leq \frac{80}{3}\epsilon
\end{eqnarray*}
for all $x\in X$.

{\bf Proof.} Consider $A$ to be ${\bf C}$ in Theorem 2.$\Box$

The next result is a generalization of the main theorem of $\cite{G-S}$.

{\bf Theorem 3.} Suppose $(X,\perp)$ is an orthogonality space and $(Y, \|.\|)$ is a real Banach space. Let $f_1, f_2, f_3:X\to Y$ be mappings fulfilling 
\begin{equation}
\|f_1(x+y)-f_2(x)-f_3(y)\|\leq\epsilon
\end{equation}
for some $\epsilon$ and for all $x,y\in X$ with $x\perp y$. Then there exists exactly a quadratic mapping $Q:X\to Y$ and an additive mapping $T:X\to Y$ such that
\begin{eqnarray*}
\|f_1(x)-f_1(0)-Q(x)-T(x)\|\leq \frac{68}{3}\epsilon\\
\|f_2(x)-f_2(0)-Q(x)-T(x)\|\leq \frac{80}{3}\epsilon\\
\|f_3(x)-f_3(0)-Q(x)-T(x)\|\leq \frac{80}{3}\epsilon
\end{eqnarray*}
for all $x\in X$.

{\bf Proof.} Use the same reasoning as in the proof of Theorem 1 with $a=b=1$ and applying Corollaries 1 and 2.$\Box$

\end{document}